\documentclass[a4paper,12pt,final]{amsart}
\usepackage{times,a4wide,mathrsfs,amssymb,amsmath,amsthm,enumerate,xypic,tikzsymbols,dsfont}

\newcommand{\C}{\mathbb{C}}
\newcommand{\ZZ}{\mathbb{Z}}

\newcommand{\QQ}{\mathbb{Q}}
\newcommand{\NN}{\mathbb{N}}
\newcommand{\PP}{\mathbb{P}}

\newcommand{\A}{\mathbb{A}}
\newcommand{\OO}{\mathcal O}
\newcommand{\Ss}{\mathcal S}

\newcommand{\Sy}{\mathfrak S}

\newcommand{\XX}{\mathcal X}
\newcommand{\YY}{\mathcal Y}
\newcommand{\UU}{\mathcal U}

\newcommand{\AAA}{\mathcal A}

\newcommand{\TT}{\mathcal T}
\newcommand{\CC}{\mathcal C}
\newcommand{\Zz}{\mathcal Z}

\newcommand{\MM}{\mathcal M}
\newcommand{\PPP}{\mathcal P}
\newcommand{\Rr}{\mathcal R}

\newcommand{\FF}{\mathcal F}

\newcommand{\pic}{\hbox{Pic}}

\newcommand{\rom}{\romannumeral}

\newcommand{\one}{\mathds{1}}

\DeclareMathOperator{\Pic}{Pic}
\DeclareMathOperator{\ima}{Im}

\DeclareMathOperator{\sym}{Sym}
\DeclareMathOperator{\Gr}{Gr}
\DeclareMathOperator{\dec}{Dec}

\newtheorem{theorem}{Theorem}[section]
\newtheorem{claim}[theorem]{Claim}

\newtheorem{lemma}[theorem]{Lemma}
\newtheorem{sublemma}[theorem]{Sublemma}
\newtheorem{corollary}[theorem]{Corollary}
\newtheorem{proposition}[theorem]{Proposition}

\newtheorem{remark}[theorem]{Remark}
\newtheorem{definition}[theorem]{Definition}
\newtheorem{convention}{Conventions}
\newtheorem{question}[theorem]{Question}
\newtheorem{notation}[theorem]{Notation}

\newtheorem{nonumbering}{Theorem}

\newtheorem{nonumberingc}{Corollary}

\newtheorem{nonumberingt}{Acknowledgements}

\begin{document}

\author[Robert Laterveer]
{Robert Laterveer}

\address{Institut de Recherche Math\'ematique Avanc\'ee,
CNRS -- Universit\'e 
de Strasbourg,\
7 Rue Ren\'e Des\-car\-tes, 67084 Strasbourg CEDEX,
FRANCE.}
\email{robert.laterveer@math.unistra.fr}

\title{Algebraic cycles and Fano threefolds of genus 10}

\begin{abstract} We show that prime Fano threefolds $Y$ of genus 10 have a multiplicative Chow--K\"unneth decomposition, in the sense of Shen--Vial. As a 
consequence, a certain tautological subring of the Chow ring of powers of $Y$ injects into cohomology.
 \end{abstract}

\thanks{Supported by ANR grant ANR-20-CE40-0023.}


\maketitle

\section{Introduction}

Given a smooth projective variety $Y$ over $\C$, let 
  \[ A^i(Y):=CH^i(Y)_{\QQ}\] 
 denote the Chow groups of $Y$ (i.e. the groups of codimension $i$ algebraic cycles on $Y$ with $\QQ$-coefficients, modulo rational equivalence). The intersection product defines a ring structure on $A^\ast(Y)=\bigoplus_i A^i(Y)$, the Chow ring of $Y$ \cite{F}. In the case of K3 surfaces, this ring has a remarkable property:

\begin{theorem}[Beauville--Voisin \cite{BV}]\label{bv} Let $S$ be a K3 surface. 
The $\QQ$-subalgebra
  \[   \bigl\langle  A^1(S), c_j(S) \bigr\rangle\ \ \ \subset\ A^\ast(S) \]
  injects into cohomology under the cycle class map.
  \end{theorem}

The Chow ring of abelian varieties also exhibits particular behaviour: there is a multiplicative splitting \cite{Beau}.
Motivated by the cases of K3 surfaces and abelian varieties, Beauville \cite{Beau3} has conjectured that for certain special varieties, the Chow ring should admit a multiplicative splitting (and a certain subring should inject into cohomology). To make concrete sense of Beauville's elusive ``splitting property conjecture'', Shen--Vial \cite{SV} have introduced the concept of {\em multiplicative Chow--K\"unneth decomposition\/}; 
we will abbreviate this to ``MCK decomposition'' (for the precise definition, cf. section \ref{ss:mck} below).

It is something of a challenge to understand precisely which varieties admit an MCK decomposition. To give an idea of what is known: hyperelliptic curves have an MCK decomposition \cite[Example 8.16]{SV}, but the very general curve of genus $\ge 3$ does not have an MCK decomposition \cite[Example 2.3]{FLV2}; K3 surfaces have an MCK decomposition, but certain high degree surfaces in $\PP^3$ do not have an MCK decomposition (cf. the examples given in \cite{OG}). In this note, we will focus on Fano threefolds and ask the following question:

\begin{question}\label{ques} Let $X$ be a Fano threefold with Picard number 1. Does $X$ admit an MCK decomposition ?
\end{question} 

The restriction on the Picard number is necessary to rule out a counterexample of Beauville \cite[Examples 9.1.5]{Beau3}. The answer to Question \ref{ques} is affirmative for cubic threefolds \cite{Diaz}, \cite{FLV2}, for intersections of 2 quadrics \cite{2q}, for intersections of a quadric and a cubic \cite{55} and for prime Fano threefolds of genus 8 \cite{g8}.

The main result of this note answers Question \ref{ques} for one more family:

\begin{nonumbering}[=Theorem \ref{main}] Let $Y$ be a prime Fano threefold of genus 10. Then $Y$ has a multiplicative Chow--K\"unneth decomposition.
\end{nonumbering}

The argument proving Theorem \ref{main} is based on the connections between $Y$ and a certain genus 2 curve, and between $Y$ and an index 2 Fano threefold $Z$ (cf. Theorem \ref{FY}). The work of Kuznetsov \cite{Ku062}, \cite{Ku2}, \cite{KPS}, building these connections on a categorical level inside the set-up of {\em homological projective duality\/}, allows to establish the instances of the {\em Franchetta property\/} that are needed to prove the theorem.

Reaping the fruits of Theorem \ref{main}, we obtain a result concerning the {\em tautological ring\/}, which is a certain subring of the Chow ring of powers of $Y$:

\begin{nonumberingc}[=Corollary \ref{cor1}] Let $Y$ be a prime Fano threefold of genus 10, and $m\in\NN$. Let
  \[ R^\ast(Y^m):=\bigl\langle (p_i)^\ast(h), (p_{ij})^\ast(\Delta_Y)\bigr\rangle\ \subset\ \ \ A^\ast(Y^m)   \]
  be the $\QQ$-subalgebra generated by pullbacks of the polarization $h\in A^1(Y)$ and pullbacks of the diagonal $\Delta_Y\in A^3(Y\times Y)$. 
  The cycle class map induces injections
   \[ R^\ast(Y^m)\ \hookrightarrow\ H^\ast(Y^m,\QQ)\ \ \ \hbox{for\ all\ }m\in\NN\ .\]
   \end{nonumberingc}

 This is the kind of injectivity result that motivated Beauville's work on the ``splitting property conjecture'' \cite{Beau3}.
 To paraphrase Corollary \ref{cor1},  one could say that genus 10 Fano threefolds behave like hyperelliptic curves from the point of view of intersection theory (cf. Remark \ref{tava} below).

 \vskip0.6cm

\begin{convention} In this article, the word {\sl variety\/} will refer to a reduced irreducible scheme of finite type over $\C$. A {\sl subvariety\/} is a (possibly reducible) reduced subscheme which is equidimensional. 

{\bf All Chow groups will be with rational coefficients}: we will denote by $A_j(Y)$ the Chow group of $j$-dimensional cycles on $Y$ with $\QQ$-coefficients; for $Y$ smooth of dimension $n$ the notations $A_j(Y)$ and $A^{n-j}(Y)$ are used interchangeably. 
The notation $A^j_{hom}(Y)$ will be used to indicate the subgroup of homologically trivial cycles.
For a morphism $f\colon X\to Y$, we will write $\Gamma_f\in A_\ast(X\times Y)$ for the graph of $f$.

The contravariant category of Chow motives (i.e., pure motives with respect to rational equivalence as in \cite{Sc}, \cite{MNP}) will be denoted 
$\MM_{\rm rat}$.
\end{convention}

\section{Prime Fano threefolds of genus 10} The classification of Fano threefolds is one of the glories of twentieth century algebraic geometry \cite{IP}. Fano threefolds that are {\em prime\/} (i.e. with Picard group of rank $1$ generated by the canonical divisor) come in 10 explicitly described families. In this paper we will be concerned with one of these families:

\begin{theorem}[Mukai \cite{Mu1}]\label{muk} Let $Y$ be a prime Fano threefold (i.e., a smooth projective Fano threefold with $\Pic(Y)=\ZZ[K_Y]$), of genus $10$. Then $Y$ is a dimensionally transverse intersection 
  \[  Y= G_2 \Gr(2,7)\cap \PP^{11}\ \ \ \subset\ \PP^{13}\ ,\]
  where $G_2 \Gr(2,7)$ is the minimal compact homogeneous space for the simple algebraic group of type $G_2$ (the variety $G_2 \Gr(2,7)$ can be realized as the zero locus of a section of a certain vector bundle on the Grassmannian $\Gr(2,7)$).
  
Conversely, any smooth dimensionally transverse intersection $G_2 \Gr(2,7)\cap \PP^{11}$ is a prime Fano threefold of genus $10$.

The Hodge diamond of $Y$ is
  \[ \begin{array}[c]{ccccccc}
      &&&1&&&\\
       &&0&&0&&\\
             &0&&1&&0&\\
                    0&&2&&2&&0\\
      &0&&1&&0&\\
      &&0&&0&&\\
      &&&1&&&\\
      \end{array}\]
  \end{theorem}

\begin{proof} The ``conversely'' statement is just because $G_2 \Gr(2,7)  $ is a Fano variety of dimension $5$ with Picard number 1, index $3$ and degree $18$; the codimension $2$ complete intersection $Y:=G_2 \Gr(2,6)\cap\PP^{11}$ thus has index $1$ and degree $d=18$ (i.e. genus $g=d/2+1=10$).
The first statement is proven in \cite{Mu1}.

To see that $h^{2,1}(Y)=2$, one can use Theorem \ref{FY}  or Theorem \ref{DY} below. 
\end{proof}

\subsection{Hilbert scheme of conics}

\begin{theorem}\label{FY} Let $Y$ be a prime Fano threefold of genus 10, and let $F:=F(Y)$ be the Hilbert scheme parametrizing conics contained in $Y$. 

\noindent
(\rom1) 
$F$ is an abelian surface, isomorphic to the intermediate Jacobian of $Y$.

\noindent
(\rom2) there is a genus 2 curve $C$, geometrically associated to $Y$, such that $F$ is isomorphic to the Jacobian of $C$.

\noindent
(\rom3) There exists $P\in A^2(Y\times F)$ inducing an isomorphism
  \[ P_\ast\colon\ \ H^3(Y,\QQ)\ \xrightarrow{\cong}\ H^1(F,\QQ)\ .\]
  
 \noindent
 (\rom4) There exists a Fano threefold $Z$ of Picard number 1, index 2 and degree 4 (i.e., $Z$ is a complete intersection of 2 quadrics in $\PP^5$), such that the Fano surface $F_1(Z)$ of lines in $Z$ is isomorphic to $F$:
   \[ F_1(Z)\cong F(Y)\ .\]
  \end{theorem}

\begin{proof} Item (\rom1) is proven for general $Y$ in \cite[Proposition 3]{IM}; the extension to arbitrary $Y$ is done in \cite[Theorem 1.1.1]{KPS}.

\noindent
(\rom2) This follows from the proof of \cite[Proposition 3]{IM} in case $Y$ is general; for arbitrary $Y$ it is part of the proof of \cite[Proposition B.5.1]{KPS}. Recall that $Y$ is a codimension 2 linear section of the homogeneous variety $G_2 \Gr(2,7)$. The pencil of hyperplanes containing $Y$ has 6 singular elements (because the projective dual to $G_2 \Gr(2,7)$ is a sextic hypersurface). The curve $C$ is constructed as the double cover $C\to\PP^1$ branched in these 6 points.

\noindent 
(\rom3) This follows from the Abel--Jacobi isomorphism (\rom1); alternatively one could use Theorem \ref{DY}.

\noindent
(\rom4) This is part of a more general (and rather mysterious) phenomenon linking certain index 1 Fano threefolds and certain index 2 Fano threefolds, that was first discovered by Kuznetsov \cite[Section 4.2]{Ku2}, \cite[Proposition B.5.1]{KPS}.

There is not a unique index 2 Fano threefold $Z$ associated to $Y$ (this is only true up to projective isomorphism), but there is a canonical way of giving a $Z$ as in (\rom4), as explained in \cite[Section B.5]{KPS}: let $\lambda_0,\ldots, \lambda_5$ denote the branch points of the double cover $C\to\PP^1$. Choosing an embedding of $\A^1$ in $\PP^1$ such that the six points $\lambda_j$ are contained in $\A^1$, let us write $\lambda_j\in\A^1$ for the affine coordinates of these six points. The intersection of quadrics
  \begin{equation}\label{lambda}  Z:= \Bigl\{  [x_0,x_1,\ldots,x_5]\in\PP^5\,\Big\vert\,  x_0^2+\cdots+x_5^2= \lambda_0 x_0^2+\cdots +\lambda_5 x_5^2=0\Bigr\}\ \ \ \subset\ \PP^5 \end{equation}
  is smooth (and hence it is a Fano threefold as in (\rom4)), and the genus 2 curve $C_Z$ (naturally associated to $Z$ by looking at the 6 singular quadrics in the pencil defining $Z$) is isomorphic to $C$. Since it is known that $F_1(Z)\cong\hbox{Jac}(C_Z)$ \cite{Reid}, this gives the required isomorphism $ F_1(Z)\cong F(Y)$.
\end{proof}

\begin{remark}\label{!} Although they are related, there is an important difference between the index 2 Fano threefolds $Z$ (as in Theorem \ref{FY}(\rom4)) and the index 1 Fano threefolds $Y$.
For the threefolds $Z$ there is a Torelli theorem, i.e. $Z$ is uniquely determined by $F_1(Z)$ \cite[Section 3.6]{Deb}.
On the other hand, the (generic) Torelli theorem fails for prime Fano threefolds of genus 10: the threefold $Y$ is {\em not\/} determined by the surface $F(Y)$. Indeed, the moduli space $\MM_{10}$ (of genus 10 prime Fano threefolds) has dimension 10, whereas
the moduli space $\MM_2$ of genus 2 curves has dimension 3, and so the surjective morphism
  \[ \MM_{10}\ \to\ \MM_2 \]
  (sending a Fano threefold $Y$ to the curve $C$) has generic fiber of dimension 7. 

In this context, Iliev and Manivel have proven a certain ``modified Torelli statement'': the genus 10 prime Fano threefold $Y$ containing a fixed K3 surface $S$ as hyperplane section is uniquely determined by the image of $F(Y)$ in the Hilbert scheme $S^{[2]}$ \cite[Theorem 6]{IM}.
\end{remark}

\subsection{The $Y$-$F(Y)$ relation}

\begin{proposition}\label{YF} Let $Y$ be a prime Fano threefold of genus 10, and let $F:=F(Y)$ be the Hilbert scheme of conics in $Y$. There is an isomorphism of Chow motives
    \[   h(Y^{(2)})\cong h(F)(-2)\oplus h(Y) \oplus h(Y)(-3)\ \ \ \hbox{in}\ \MM_{\rm rat}\ .\]
\end{proposition}

\begin{proof} (NB: the symmetric product $Y^{(2)}$ is not smooth, but it is a projective Alexander scheme in the sense of \cite{KimA} and so it makes sense to speak about the motive of $Y^{(2)}$ in the category $\MM_{\rm rat}$. If one prefers, one may think of the motive $h(Y^{(2)})$ as $(h(Y)\otimes h(Y))^{\Sy_2}$.)

Letting $h\in A^1(Y)$ denote a hyperplane section (with respect to the embedding $Y\subset\PP^{13}$ given by Theorem \ref{muk}), let us write
    \begin{equation}\label{ck}  \begin{split}    \pi^0_Y&:= {1\over 18}\, h^3\times Y\ ,\\
                              \pi^2_Y&:= {1\over 18}\, h^2\times h\ ,\\
                               \pi^4_Y&:= {1\over 18}\, h\times h^2\ ,\\
                                \pi^6_Y&:= {1\over 18}\, Y\times h^3\ ,\\      
                                \pi^3_Y&:= \Delta_Y-\sum_{j\not=3} \pi^j_Y\ \ \ \in\ A^3(Y\times Y)\ ,\\
                                \end{split}\end{equation}
                                and $h^j(Y):=(Y,\pi^j_Y,0)\in\MM_{\rm rat}$.
                                (This is the CK decomposition which we will prove to be MCK in Theorem \ref{main}.)
                                
  As we have seen (Theorem \ref{FY}), there is an isomorphism
    \[ P_\ast\colon\ \ H^3(Y,\QQ)\ \xrightarrow{\cong}\ H^1(F,\QQ)\ .\]
    Since both $Y$ and $F$ verify the Lefschetz standard conjecture, the inverse isomorphism is also induced by a correspondence. (This is well-known, cf. for instance \cite[Proof of Proposition 1.1]{V4}, where I first learned this.) It follows that there is an isomorphism of homological motives
    \[ P\colon\ \ h^3(Y)\    \xrightarrow{\cong}\ h^1(F)(-1)         \ \ \ \hbox{in}\ \MM_{\rm hom}\ .\]
 Since $Y$ and $F$ are Kimura finite-dimensional (in the sense of \cite{Kim}), this can be upgraded to an isomorphism of Chow motives
  \begin{equation}\label{isoYF}  P\colon\ \ h^3(Y)\    \xrightarrow{\cong}\ h^1(F)(-1)         \ \ \ \hbox{in}\ \MM_{\rm rat}\ .\end{equation}
Writing $h(Y)\cong h^1(F)(-1)\oplus \one\oplus \one(-1)\oplus \one(-2)\oplus \one(-3)$ and taking the symmetric power, one obtains
  \begin{equation}\label{thisiso} \begin{split}  h(Y^{(2)})\cong  \sym^2 h^1(F)(-2)\oplus h^1(F)(-1) \oplus h^1(F)(-2) \oplus h^1(F)(-3)\oplus h^1(F)(-4)&\\ \oplus \bigoplus \one(\ast)\ \ \ \hbox{in}\ &\MM_{\rm rat}\ .\\
  \end{split} \end{equation}

 On the other hand, $F$ being an abelian surface its motive decomposes
    \begin{equation}\label{andthis} h(F)= \one \oplus h^1(F) \oplus \sym^2 h^1(F) \oplus h^1(F)(-1) \oplus \one(-2)\ \ \ \hbox{in}\ \MM_{\rm rat}\ \end{equation}
    (cf. for instance \cite[Section 5]{Sc}).
    Combining \eqref{thisiso} and \eqref{andthis}, we obtain the isomorphism of the proposition.
    \end{proof}

 \begin{remark} Proposition \ref{YF} is formally similar to the relation between a cubic hypersurface $Y\subset\PP^{n+1}$ and its Fano variety of lines $F:=F(Y)$: in this case, it is known that
   \[  h(Y^{(2)})\cong h(F)(-2)\oplus h(Y)\oplus h(Y)(-n)\ \ \ \hbox{in}\ \MM_{\rm rat}\ \]
   \cite[Theorem 3]{FLV3}. This is the motivic version of the $Y$-$F(Y)$ relation in the Grothendieck ring of varieties that was discovered by Galkin--Shinder \cite{GS}. 
 \end{remark}

\subsection{The $Z$-$F(Z)$ relation}

\begin{proposition}\label{ZF} Let $Z\subset\PP^5$ be a smooth complete intersection of 2 quadrics, and let $F:=F_1(Z)$ be the Hilbert scheme of lines in $Z$. There is an isomorphism of Chow motives
    \[   h(Z^{(2)})\cong h(F)(-2)\oplus h(Z) \oplus h(Z)(-3)\ \ \ \hbox{in}\ \MM_{\rm rat}\ .\]
\end{proposition}

\begin{proof} 
The argument is the same as that of Proposition \ref{YF}, given that $F$ is an abelian surface and that the universal line $P\subset Z\times F$ induces an Abel--Jacobi isomorphism
  \[   P_\ast\colon\ \ H^3(Z,\QQ)\ \xrightarrow{\cong}\ H^1(F,\QQ) \]
 \cite[Theorem 4.14]{Reid}. 
  \end{proof}

\subsection{Conics on $G$} Let us write $G:= G_2 \Gr(2,7)$ for the $G_2$-Grassmannian. We will need to understand conics on $G$:

\begin{proposition}\label{FG}
 Let $G:= G_2 \Gr(2,7)$, and let $F(G)$ be the Hilbert scheme of conics (with respect to the Pl\"ucker embedding) contained in $G$.
The scheme $F(G)$ is a smooth projective spherical variety of dimension 8. In particular, $A^\ast_{hom}(F(G))=0$.
\end{proposition}

\begin{proof} The variety $F(G)$ is isomorphic to the so-called {\em Cayley Grassmannian\/} (cf. \cite[Section 7.1]{Man}, where this observation is attributed to Kuznetsov). The Cayley Grassmannian $CG$ is known to be a smooth projective spherical variety \cite[Theorem 1.1]{Man0}, hence it has trivial Chow groups. (Actually, it is known that the Chow groups of $CG$ 
are isomorphic to those of the Grassmannian $\Gr(2,6)$ \cite[Theorem 1.2]{Man0}.)
\end{proof}

\subsection{Derived category}

 \begin{theorem}\label{DY} Let $Y$ be a prime Fano threefold of genus 10. The derived category of $Y$ admits a semi-orthogonal decomposition
   \[ D^b(Y)=    \bigl\langle D^b(C), \OO_Y,\UU^\ast\bigr\rangle   \ ,\] 
   where $C$ is the genus 2 curve of Theorem \ref{YF} and $\UU$ is the restriction of the tautological rank 2 subbundle on $\Gr(2,7)$.
 \end{theorem}
 
 \begin{proof} This is proven in \cite[Section 6.4]{Ku062}, as an instance of {\em homological projective duality\/}. It will be important for us to understand how the curve $C$ and the semi-orthogonal decomposition are constructed. As in loc. cit., let $V$ be the 14-dimensional vector space such that $G:=G_2 \Gr(2,7)$ is embedded in $\PP(V)$ (writing $\Gr(2,W)=\Gr(2,7)$, this $V$ arises from the decomposition $\wedge^2 V = W^\ast\oplus V$ into $G_2$-representations). The projective dual $G^\vee\subset\PP(V^\ast)$ of $G$ is a sextic hypersurface, with singular locus $G^\vee_{\rm sing}$ of dimension 10. Let 
   \[ h\colon \ \ H\ \to\ \PP(V^\ast)\setminus G^\vee_{\rm sing}\] 
 be the double cover ramified along $G^\vee\setminus G^\vee_{\rm sing}$.
  As explained in loc. cit., there is a certain sheaf of algebras $\AAA_H$ on $H$ such that the non-commutative variety
   \[   (H,\AAA_H)   \]
   is HPD dual to $G_2 \Gr(2,7)$. As proven in loc. cit., this entails in particular that given an $r$-dimensional linear subspace $L\subset V^\ast$, one can relate the derived categories of the linear sections $H_r:= h^{-1}(\PP(L))$ and $G_r:= G\cap \PP(L^\ast)$. Taking an $L$ of dimension $r=2$ and such that $H_2$ (and hence $G_2$) is smooth and dimensionally transverse, one has that $H_2$ is a genus 2 curve, and $G_2$ is a prime Fano threefold of genus 10 (and every prime Fano threefold of genus 10 arises in this way, cf. Theorem \ref{muk}). This gives the semi-orthogonal decomposition as indicated.
    
   (For later use, we further observe that taking a linear subspace $L$ of dimension $r=3$ and intersecting smoothly and transversely, the varieties $G_3$ and $H_3$ are K3 surfaces of degree 18 resp. 2 that are twisted derived equivalent: as shown in loc. cit., one has
     \begin{equation}\label{K3} D^b(G_3)\cong D^b(H_3,\alpha)\ ,\end{equation}
     where $\alpha$ is a Brauer class.
     Moreover, the general K3 surfaces of degree 18 and of degree 2 arise in this way.)
       \end{proof}

\section{MCK decomposition}
\label{ss:mck}

\begin{definition}[Murre \cite{Mur}] Let $X$ be a smooth projective variety of dimension $n$. We say that $X$ has a {\em CK decomposition\/} if there exists a decomposition of the diagonal
   \[ \Delta_X= \pi^0_X+ \pi^1_X+\cdots +\pi_X^{2n}\ \ \ \hbox{in}\ A^n(X\times X)\ ,\]
  such that the $\pi^i_X$ are mutually orthogonal idempotents and $(\pi_X^i)_\ast H^\ast(X,\QQ)= H^i(X,\QQ)$.
  
  (NB: ``CK decomposition'' is shorthand for ``Chow--K\"unneth decomposition''.)
\end{definition}

\begin{remark} Murre has conjectured that any smooth projective variety should have a CK decomposition \cite{Mur}, \cite{J4}. 
\end{remark}

\begin{definition}[Shen--Vial \cite{SV}] Let $X$ be a smooth projective variety of dimension $n$, and let $\Delta_X^{sm}\in A^{2n}(X\times X\times X)$ denote the class of the small diagonal
  \[ \Delta_X^{sm}:=\bigl\{ (x,x,x)\ \vert\ x\in X\bigr\}\ \subset\ X\times X\times X\ .\]
  An {\em MCK decomposition\/} is defined as a CK decomposition $\{\pi_X^i\}$ of $X$ that is {\em multiplicative\/}, i.e. it satisfies
  \[ \pi_X^k\circ \Delta_X^{sm}\circ (\pi_X^i\times \pi_X^j)=0\ \ \ \hbox{in}\ A^{2n}(X\times X\times X)\ \ \ \hbox{for\ all\ }i+j\not=k\ .\]
  
 (NB: ``MCK decomposition'' is shorthand for ``multiplicative Chow--K\"unneth decomposition''.) 
  \end{definition}
  
  \begin{remark} The small diagonal (when considered as a correspondence from $X\times X$ to $X$) induces the {\em multiplication morphism\/}
    \[ \Delta_X^{sm}\colon\ \  h(X)\otimes h(X)\ \to\ h(X)\ \ \ \hbox{in}\ \MM_{\rm rat}\ .\]
 Let us assume $X$ has a CK decomposition
  \[ h(X)=\bigoplus_{i=0}^{2n} h^i(X)\ \ \ \hbox{in}\ \MM_{\rm rat}\ .\]
  By definition, this decomposition is multiplicative if for any $i,j$ the composition
  \[ h^i(X)\otimes h^j(X)\ \to\ h(X)\otimes h(X)\ \xrightarrow{\Delta_X^{sm}}\ h(X)\ \ \ \hbox{in}\ \MM_{\rm rat}\]
  factors through $h^{i+j}(X)$.
  
  If $X$ has an MCK decomposition, then setting
    \[ A^i_{(j)}(X):= (\pi_X^{2i-j})_\ast A^i(X) \ ,\]
    one obtains a bigraded ring structure on the Chow ring: that is, the intersection product sends $A^i_{(j)}(X)\otimes A^{i^\prime}_{(j^\prime)}(X) $ to  $A^{i+i^\prime}_{(j+j^\prime)}(X)$.
    
      It is conjectured that for any $X$ with an MCK decomposition, one has
    \[ A^i_{(j)}(X)\stackrel{??}{=}0\ \ \ \hbox{for}\ j<0\ ,\ \ \ A^i_{(0)}(X)\cap A^i_{hom}(X)\stackrel{??}{=}0\ ;\]
    this is related to Murre's conjectures B and D, that have been formulated for any CK decomposition \cite{Mur}. In particular, this would imply that the subring $A^\ast_{(0)}(X)$ injects into cohomology under the cycle class map.
    \end{remark}

    \begin {remark}
 The property of having an MCK decomposition is motivated by, and closely related to, Beauville's ``splitting property' conjecture'' \cite{Beau3}. 
  To give an idea of what is known: hyperelliptic curves have an MCK decomposition \cite[Example 8.16]{SV}, but the very general curve of genus $\ge 3$ does not have an MCK decomposition \cite[Example 2.3]{FLV2}. It has been conjectured that all hyperk\"ahler varieties have an MCK decomposition. For a more thorough discussion, and for more examples of varieties with an MCK decomposition, we refer to \cite[Section 8]{SV}, as well as \cite{V6}, \cite{SV2}, \cite{FTV}, \cite{37}, \cite{38}, \cite{39}, \cite{40}, \cite{44}, \cite{FLV2}, \cite{60}, \cite{55}, \cite{59}.
    \end{remark}

 \section{Franchetta property}
 
 \begin{definition} Let $\XX\to B$ be a smooth projective morphism, where $\XX, B$ are smooth quasi-projective varieties. We say that $\XX\to B$ has the {\em Franchetta property in codimension $j$\/} if the following holds: for every $\Gamma\in A^j(\XX)$ such that the restriction $\Gamma\vert_{X_b}$ is homologically trivial for the very general $b\in B$, the restriction $\Gamma\vert_b$ is zero in $A^j(X_b)$ for all $b\in B$.
 
 We say that $\XX\to B$ has the {\em Franchetta property\/} if $\XX\to B$ has the Franchetta property in codimension $j$ for all $j$.
 \end{definition}
 
 This property is studied in \cite{PSY}, \cite{BL}, \cite{FLV}, \cite{FLV3}.
 
  \begin{definition} Let $\XX\to B$ be a family as above, with $X:=X_b$ a fiber. We will write
   \[ GDA^j_B(X):=\ima\Bigl( A^j(\XX)\to A^j(X)\Bigr) \]
   for the subgroup of {\em generically defined cycles}. 
   (In a context where it is clear to which family we are referring, the index $B$ will sometimes be dropped from the notation.)
  \end{definition}
  
  With this definition, the Franchetta property amounts to saying that $GDA^\ast(X)$ injects into cohomology, under the cycle class map.

 \subsection{Franchetta property for $Y$}
 
 \begin{notation}\label{not} Let $G$ be the $G_2$-Grassmannian $G:=G_2 \Gr(2,7)$, and let $\OO_G(1)$ be the polarization corresponding to the Pl\"ucker embedding
 $G\subset\PP^{13}$. Let 
   \[ B\ \subset\ \bar{B}:=\PP H^0(G,\OO_G(1)^{\oplus 2})\] 
   denote the Zariski open subset parametrizing smooth dimensionally transverse complete intersections, and let 
   \[ \YY\ \to\ B \]
   denote the universal family of smooth $3$-dimensional complete intersections (in view of Theorem \ref{muk}, this is the universal family of prime Fano threefolds of genus $10$).
  \end{notation}
 
 \begin{proposition}\label{Fr} Let $\YY\to B$ be the universal family of prime Fano threefolds of genus $10$ (Notation \ref{not}). 
 The family $\YY\to B$ has the Franchetta property.
  \end{proposition}
 
 \begin{proof} We give two different proofs of this proposition. For the first proof,
 let $\bar{\YY}\subset\bar{B}\times G$ denote the projective closure. As the line bundle $\OO_G(1)$ is base point free, the projection $\bar{\YY}\to G$ is a $\PP^r$-fibration.
 Using the projective bundle formula, it is readily checked (cf. for instance \cite[Proof of Lemma 1.1]{PSY}) that
   \[ GDA^\ast_B(Y) =\ima\bigl( A^\ast(G)\to A^\ast(Y)\bigr)\ .\]
   Since $A^j_{hom}(Y)=0$ for $j\not=2$, it only remains to ascertain that the cycle class map induces injections
     \begin{equation}\label{in} \ima\bigl(  A^\ast(G)\to A^\ast(Y)\bigr)\ \to\ H^\ast(Y,\QQ)\ .\end{equation}
  But the Chow ring of the $G_2$-Grassmannian $G$ is as small as can be for a smooth projective variety: indeed, the fivefold $G$ admits a full exceptional collection with 6 exceptional objects \cite[Section 6.4]{Ku062}, which means that $G$ is a {\em minifold\/} in the sense of \cite{GK+}, \cite[Section 1.4]{Ku3}.
  Taking Chow groups, we find that $A^\ast(G)=\QQ^6$ and so $A^j(G)=\QQ[h^j]$ for all $j$. This implies that \eqref{in} is injective, and ends the first proof.

   For the second proof of the proposition, we explore the relation with the genus 2 curve $C$ given by the HPD framework
   of Theorem \ref{DY}. This theorem gives an isomorphism
   \[ R_\ast\colon\ \ A^2(Y)\ \xrightarrow{\cong}\ A^1(C)\ .\]
   Even better, this isomorphism exists {\em universally\/}: let us write $\CC\to B$ for the universal family of smooth one-dimensional linear sections of $H$ (where $H$ is the double cover defined in the proof of Theorem \ref{DY}). The correspondence $R$ is generically defined (i.e. there exists $\Rr\in A^2(\YY\times_B \CC)$ such that the fiberwise restriction $\Rr\vert_b$ is the correspondence $R\in A^2(Y\times C)$); this is because the semi-orthogonal decomposition of Theorem \ref{DY} exists universally (cf.    \cite[Proof of Theorem 1.2]{Ku062}).
   It follows that $R_\ast$ sends generically defined cycles to generically defined cycles, and so there is a commutative diagram
   \[ \begin{array}[c]{ccc}
                     GDA^2_B(Y)  & \xrightarrow{R_\ast} & GDA^1_B(C)\\
                     &&\\
                     \downarrow&&\downarrow\\
                     &&\\
                  \QQ\cong   H^4(Y,\QQ)& \xrightarrow{R_\ast} & H^2(C,\QQ)\cong\QQ\\   
                  \end{array}\]
                  (where vertical arrows are cycle class maps).
                  
      To prove the proposition, it thus suffices to prove the Franchetta property for the family of curves $\CC\to B$. Let $\bar{\CC}\subset \bar{B}\times\bar{H}$ denote the projective closure, then $\bar{\CC}\to \bar{H}$ is a $\PP^s$-fibration, and so (as above) the projective bundle formula (plus the fact that $C$ is contained in the smooth quasi-projective variety $H$, by construction) gives
      \[  GDA^1_B(C)=\ima\bigl( A^1(H)\to A^1(C)\bigr)\ .\]
      Let us check that the right-hand side is one-dimensional. In view of the usual spread lemma (cf. \cite[Lemma 3.2]{Vo}), it suffices to prove this for the very general $C$, and so we may suppose $C$ is contained in a 2-dimensional smooth linear section $C:=H_2 \subset H_3\subset H$. Such a surface $H_3$ is a degree 2 K3 surface. What's more, in view of Lemma \ref{pic} below, we may suppose $H_3$ has Picard number 1, and so 
      \[ \ima \bigl( A^1(H)\to A^1(C)\bigr)  = \ima \bigl( A^1(H_3)\to A^1(C)\bigr) =\QQ[K_C]\ .\]
      The proposition is now proven, modulo the following lemma:
      
      \begin{lemma}\label{pic} The very general curve $C\subset H$ is contained in a linear section $H_3\subset H$, where $H_3$ is a K3 surface with $\pic(H_3)=\ZZ$.
      \end{lemma}
      
      To prove the lemma, we return to the HPD framework of Theorem \ref{DY}: thanks to Mukai, we know that the very general K3 surface of degree 18 arises as a linear section $G_3$ of the $G_2$-Grassmannian $G$, and so the very general Fano threefold $Y$ has an anticanonical section $G_3\subset Y$ with $\pic(G_3)\cong\ZZ$. On the HPD dual side, the genus 2 curve $C\subset H$ associated to $Y$ is contained in the K3 surface $H_3\subset H$ which is twisted derived equivalent to $G_3$, cf. equality \eqref{K3}. Twisted derived equivalent K3 surfaces are isogenous (and even have isomorphic Chow motives, cf. \cite{Huy}) and so $\pic(H_3)\cong\ZZ$.
                  \end{proof}

 \subsection{Franchetta property for $Y\times Y$}

 \begin{proposition}\label{Fr2} Let $\YY\to B$ be as in Notation \ref{not}. The family $\YY\times_B \YY\to B$ has the Franchetta property.
  \end{proposition} 
  
  \begin{proof} To prove this, we move to the HPD dual side, as in the second proof of Proposition \ref{Fr}.
Kuznetsov's work (Theorem \ref{DY}) gives an isomorphism of motives
  \[ h^3(Y)\cong h^1(C)(-1)\ \ \ \hbox{in}\ \MM_{\rm rat}\ ,\]
  where $C$ is a genus 2 curve.
 It follows that there is a
 split injection of motives
  \[ h(Y)\ \hookrightarrow\ h(C)(-1)\oplus \one \oplus \one(-1)\oplus \one(-2)\oplus \one(-3)\ \ \ \hbox{in}\ \MM_{\rm rat}\ .\]
  As we have verified in the proof of Proposition \ref{Fr}, this isomorphism and split injection are generically defined (with respect to $B$), and so one obtains in particular split injections of Chow groups
    \[ GDA^j_B(Y\times Y)\ \hookrightarrow\ GDA^{j-2}_B(C\times C)\oplus \bigoplus GDA^\ast_B(C) \oplus \QQ^t\ .\]
    Since this injection is compatible with cycle class maps, it suffices to prove the Franchetta property for the family $\CC\times_B \CC\to B$ (where $\CC\to B$ is as in the proof of Proposition \ref{Fr}). To this end, we recall that $\CC\to B$ is the universal family of smooth one-dimensional linear sections of $H$ (cf. the proof of Theorem \ref{DY}). Since $\OO_H(1)$ is very ample, this set-up verifies the property ($\ast_2$) of \cite{FLV}, and so \cite[Proposition 5.2]{FLV} implies that there is equality
    \begin{equation}\label{gda} GDA^\ast_B(C\times C) = \bigl\langle  \ima\bigl(A^\ast(H\times H)\to A^\ast(C\times C)\bigr),\, \Delta_C\bigr\rangle\ .\end{equation}
    It is left to check that the right-hand side of \eqref{gda} injects into cohomology. We need a lemma:
    
    \begin{lemma}\label{lk3} Let $\Ss\to B^\prime$ be the universal family of smooth 2-dimensional linear sections of $H$ (i.e., the fibers of $\Ss\to B^\prime$ are the degree 2 K3 surfaces $H_3$ of \eqref{K3}). The family $\Ss\times_{B^\prime}\Ss\to B^\prime$ has the Franchetta property.
     \end{lemma}
     
     \begin{proof} The HPD set-up (cf. the proof of Theorem \ref{DY}) shows that there is a ``dual'' family $\TT\to B^\prime$ over the same base, where the fibers of $\TT\to B^\prime$ are the degree 18 K3 surfaces $G_3\subset G$ of \eqref{K3}. For any $b\in B^\prime$, the fibers $S_b$ and $T_b$ are isogenous and so they have isomorphic Chow motives
     \cite{Huy}. This isomorphism of motives is universally defined (indeed, the Fourier--Mukai functor inducing the derived equivalence between $S_b$ and $T_b$ is universally defined, by construction), and so the Franchetta property for $\Ss\times_{B^\prime}\Ss\to B^\prime$ is equivalent to the Franchetta property for $\TT\times_{B^\prime}\TT\to B^\prime$. The latter property is proven in \cite[Theorem 1.5]{FLV}.     
      \end{proof}
      
      Armed with Lemma \ref{lk3}, let us further analyze the right-hand side of \eqref{gda}. Up to shrinking the base $B$ (and invoking the spread lemma \cite[Lemma 3.2]{Vo}), we may assume the curve $C\subset H$ is contained in a smooth K3 surface $S\subset H$ such that $S$ is a fiber of the family $\Ss\to B^\prime$ of Lemma \ref{lk3}.
      Then we observe that (because of the inclusions $C\subset S\subset H$) the restriction map
        \[ A^\ast(H\times H)\ \to\ A^\ast(C\times C) \]
        factors over $GDA^\ast_{B^\prime}(S\times S)$. Using \eqref{gda} and Lemma \ref{lk3}, it follows that
        \[ GDA^1_B(C\times C)= \bigl\langle (p_i)^\ast(h),\, \Delta_C\bigr\rangle \cap A^1(C\times C) \ .\]
        Since $\Delta_C$ is not decomposable in cohomology (otherwise $H^1(C,\QQ)$ would be zero), this shows the injectivity of $GDA^1_B(C\times C)\to H^2(C\times C,\QQ)$.
        Next, for the codimension 2 cycles we observe that
        \[ \begin{split} \ima\bigl( A^2(H\times H)\to A^2(C\times C)\bigr)\ \ &\subset\ \ima\bigl(   GDA_{B^\prime}^2(S\times S)\to A^2(C\times C)\bigr)\\
         &= \ \bigl\langle (p_i)^\ast(h), \, \Delta_S\vert_{C\times C}\bigr\rangle\ .\end{split}\]
        (For the last equality, we have used Lemma \ref{lk3} which gives an equality $GDA^2_{B^\prime}(S\times S)=\langle (p_i)^\ast(h),\Delta_S\rangle$.) Since $C\subset S$ is a hyperplane section, the excess intersection formula gives the equality $\Delta_S\vert_{C\times C}=\Delta_C\cdot (p_1)^\ast(h)$. Now, $h$ is proportional to the canonical divisor $K_C$ in $A^1(C)$ and so 
        \begin{equation}\label{inc} \Delta_C\cdot (p_1)^\ast(h) =   \Delta_C\cdot (p_1)^\ast(K_C) \ \ \in\ \bigl\langle (p_i)^\ast(K_C)\bigr\rangle =  \bigl\langle (p_i)^\ast(h)\bigr\rangle  \  \end{equation}  
      (this inclusion is true more generally for hyperelliptic curves \cite{Ta}, but {\em not\/} for general curves of genus $\ge 4$, cf. \cite{GG}, \cite{Yin0}).     
      It follows that 
      \[ \ima\bigl( A^2(H\times H)\to A^2(C\times C)\bigr) =\bigl\langle (p_i)^\ast(h)\bigr\rangle\ ,\]
      and so \eqref{gda} simplifies in codimension 2 to
      \[  \begin{split} GDA^2_B(C\times C) &= \bigl\langle  \ima\bigl(A^\ast(H\times H)\to A^\ast(C\times C)\bigr), \, \Delta_C\bigr\rangle    \  \cap A^2(C\times C)\\
                                                                 &=\bigl\langle (p_i)^\ast(h),\,  \Delta_C\bigr\rangle    \  \cap A^2(C\times C)\\
                                                                 &= \bigl\langle (p_i)^\ast(h) \bigr\rangle    \  \cap A^2(C\times C)\\  
                                                                 \end{split}\]
                                                  (where in the last equality we have used once more the inclusion \eqref{inc}). In view of the K\"unneth decomposition of cohomology, it is now clear that $GDA^2_B(C\times C)$ injects into cohomology. This closes the proof.              
       \end{proof}

 \begin{corollary}\label{Deltah} Let $Y$ be a genus 10 prime Fano threefold. There exist $a_j\in\QQ$ such that there is equality
   \[ \Delta_Y\cdot (p_i)^\ast(K_Y) =\sum_{j=1}^{3} a_j \, K_Y^j\times K_Y^{4-j}\ \ \ \hbox{in}\ A^4(Y\times Y)\ \ \ (i=1,2)\ .\]
  \end{corollary}     
  
  \begin{proof} One can readily find $a_j\in\QQ$ such that the equality of the corollary is true in cohomology (this is because $\Delta_Y\cdot (p_i)^\ast(K_Y) $ is the correspondence
  acting as cup product with $K_Y$; this action is non-zero only on $H^{2j}(Y,\QQ)$ which is algebraic). The Franchetta property of Proposition \ref{Fr2} then allows to lift the equality to rational equivalence.  
  \end{proof}

\begin{corollary}\label{FrF} Let $\FF\to B$ denote the universal family of Hilbert schemes of conics contained in genus 10 prime Fano threefolds. The family $\FF\to B$ has the Franchetta property.
\end{corollary}

\begin{proof} The $Y$-$F(Y)$ relation of motives of Proposition \ref{YF} is universally defined (with respect to $B$); indeed, this relation is based on the isomorphism $R\colon h^3(Y)\cong h^1(F)(-1)$ of Theorem \ref{DY}, which (as we have seen in the proof of Proposition \ref{Fr}) is universally defined. This means that there is a commutative diagram
  \[ \begin{array}[c]{ccc}
                  GDA^j_B(F) & \to & GDA^{j+2}_B(Y^{(2)}) \oplus \bigoplus GDA^\ast_B(Y) \oplus \QQ^t\\
                  &&\\
                  \downarrow&&\downarrow\\
                  &&\\
                  H^{2j}(F,\QQ) &\to& \ H^{2j+4}(Y^{(2)},\QQ)\oplus \bigoplus H^\ast(Y,\QQ) \oplus \QQ^t\ ,\\
                  \end{array}\]
                  where the horizontal arrows are injections, and the vertical arrows are the cycle class maps.
                  The right vertical arrow is injective (this follows from the Franchetta property for $Y$ and for $Y^2$, cf. Propositions \ref{Fr} and \ref{Fr2}), and so the left vertical arrow is injective as well.
\end{proof}

\subsection{Franchetta for $Z\times Z$}

\begin{notation}\label{Bcirc} Let $\YY\to B$ be the universal family as in Notation \ref{not}. Likewise, let
  \[ B_{(2,2)}\subset \PP H^0(\PP^5, \OO_{\PP^5}(2)^{\oplus 2}) \]
  be the Zariski open parametrizing smooth complete intersections of 2 quadrics, and let 
  \[    B_{(2,2)}\times \PP^5        \ \supset\ \Zz\ \to\ B_{(2,2)} \]
  denote the universal family of smooth complete intersections of 2 quadrics in $\PP^5$.
  
  The construction outlined in the proof of Theorem \ref{FY}(\rom4) gives a non-empty Zariski open $B^\circ \subset B$ and a map
   \[ B^\circ \ \to\ B_{(2,2)} \]
   associating to a prime genus 10 Fano threefold $Y$ an intersection of 2 quadrics $Z$ such that there is an isomorphism $F(Y)\cong F_1(Z)$.
   We will write
     \[ \Zz\ \to\ B^\circ \]
     for the family obtained by base change.
\end{notation}

\begin{proposition}\label{Fr2Z} The family $\Zz\times_{B^\circ}\Zz\to B^\circ$ has the Franchetta property.
\end{proposition}

\begin{proof} We will use the following motivic relation:

\begin{lemma}\label{ZY} Let $Y$ be a genus 10 Fano threefold parametrized by $B^\circ$, and let $Z$ be the associated complete intersection of 2 quadrics (cf. Theorem \ref{FY}). There is an isomorphism
of motives
  \[  \Gamma\colon\ \ h(Z)\ \xrightarrow{\cong}\ h(Y)\ \ \ \hbox{in}\ \MM_{\rm rat}\ .\]
  Moreover, the correspondence $\Gamma$ and the correspondence inducing the inverse isomorphism are generically defined with respect to $B^\circ$.
\end{lemma}

The lemma obviously implies the proposition: the correspondence $\Gamma\times\Gamma$ induces an isomorphism
  \[ GDA^\ast_{B^\circ}(Z\times Z)\ \xrightarrow{\cong}\ GDA^\ast_{B^\circ}(Y\times Y)  \]
  compatible with cycle class maps, and so Proposition \ref{Fr2Z} follows from Proposition \ref{Fr2}.
  
  To prove the lemma, as both sides are Kimura finite-dimensional it suffices to construct an isomorphism of homological motives. In addition, since $h^j(Z)\cong h^j(Y)\cong\one(\ast)$ for $j\not=3$, it suffices to construct an isomorphism between $h^3(Z)$ and $h^3(Y)$. To this end, we consider the composition
    \[   h^3(Z)\ \xrightarrow{{}^t P}\ h^1(F_1(Z))(-1)\ \xrightarrow{\cong}\ h^1(F(Y))(-1)\ \xrightarrow{Q}\ h^3(Y)\ \ \ \hbox{in}\ \MM_{\rm hom}\ ,\]
    where the middle isomorphism is induced by the isomorphism of varieties $F_1(Z)\cong F(Y)$ (Theorem \ref{YF}), the first map is defined by the transpose of the universal line (which is an isomorphism by \cite[Theorem 4.14]{Reid}), and the last map is the inverse of the isomorphism of Theorem \ref{YF}. 
   A general Hilbert schemes argument (cf. \cite[Proposition 2.11]{56} or \cite[Proposition 2.11]{cubic}) implies that $Q$ (and hence $\Gamma$) may be assumed to be generically defined; the same argument applies to the correspondence inducing the inverse isomorphism to $\Gamma$.
  \end{proof}

\begin{remark} The Franchetta property for the family  $\Zz\times_{B_{(2,2)}}\Zz\to B_{(2,2)}$ is established in \cite[Proposition 3.6(\rom2)]{2q}. However, this does not imply anything for the base-change to
$B^\circ$ (there may be more algebraic cycles on the base-changed family).

The proof of Proposition \ref{Fr2Z} relies in an essential way on the relation between $Z$ and the genus 10 Fano threefold $Y$.
\end{remark}

  \subsection{Franchetta for $Z^{(2)}\times Z$}
  
  \begin{proposition}\label{Fr3Z} Let $\Zz\to B^\circ$ be as above (Notation \ref{Bcirc}). The family
    \[   \bigl(\Zz\times_{B^\circ} \Zz\times_{B^\circ} \Zz\bigr)/\Sy_2\ \to\ B^\circ \]
    (where $\Sy_2$ acts by exchanging the first 2 factors) has the Franchetta property.
  \end{proposition}

\begin{proof}
We are going to use the Fano surface of lines $F:=F_1(Z)$. Let $\Zz\to B^\circ$ be the universal family of intersections of 2 quadrics as in Notation \ref{Bcirc}, and let
$\FF\to B^\circ$ denote the universal family of Fano surfaces of lines.
  We now make the following claim:
    
    \begin{claim}\label{cl} The family 
      \[ \FF\times_{B^\circ} \Zz\ \ \to\ B^\circ\] 
      has the Franchetta property.
    \end{claim}
    
    The claim suffices to prove Proposition \ref{Fr3Z}: indeed, Proposition \ref{ZF} gives a (generically defined) isomorphism of motives
      \[  h(Z^{(2)})\ \cong\    h(F)(-2)   \oplus h(Z) \oplus h(Z)(-3)\ \ \ \hbox{in}\ \MM_{\rm rat}\ .\]
      In particular, this gives an isomorphism of Chow groups
      \[ GDA^j_{B^\circ}(Z^{(2)}\times Z)\ \cong\ GDA_{B^\circ}^{j-2}(F\times Z)\oplus GDA^j_{B^\circ}(Z\times Z) \oplus GDA^{j-3}_{B^\circ}(Z\times Z)\ ,\]
      compatible with cycle class maps. Since we already know that $\Zz\times_{B^\circ} \Zz\to B^\circ$ has the Franchetta property (Proposition \ref{Fr2Z}), the claim thus implies Proposition \ref{Fr3Z}. 
      
      To prove Claim \ref{cl}, we borrow the argument of the closely related \cite[Proposition 3.10]{2q} (which in its turn is inspired by the work of Diaz \cite{Diaz}).
        
Let us write 
      \[ P_\PP \ \subset\   \Gr(2,6)\times \PP^5\] 
     for the universal line, and $P\subset F_1(Z)\times Z$ for the restriction of $P_\PP$ to $F_1(Z)\times Z$.

      Let $\bar{\FF}\subset\bar{B}\times \Gr(2,6)$ and $\bar{\Zz}\subset\bar{B}\times \PP^5$ denote the projective closures.
  We now consider       
the projection
    \[ \pi\colon\ \  \bar{\FF}\times_{\bar{B}} \bar{\Zz}\ \to\  \Gr(2,6)\times \PP^5\ . \]
    
As a first step towards proving Claim \ref{cl}, let us ascertain that $\pi$ has the structure of a {\em stratified projective bundle\/} (in the sense of \cite[Section 5.1]{FLV}):    
    \begin{lemma}\label{strat} 
   The morphism $\pi$
   has the structure of a $\PP^{r}$-bundle over 
     $(  \Gr(2,6)\times \PP^5\ )\setminus P_\PP$, 
     and a $\PP^{s}$-bundle over $ P_\PP$. 
       \end{lemma}

\begin{proof} Let $B_{(2,2)}$ be as in Notation \ref{Bcirc}. In \cite[Proof of Proposition 3.10]{2q}, it is proven that the family
    \[  \bar{\FF}\times_{\bar{B}_{(2,2)}} \bar{\Zz}\ \to\  \Gr(2,6)\times \PP^5\ . \]
is a $\PP^{r}$-bundle over 
     $(  \Gr(2,6)\times \PP^5\ )\setminus P_\PP$, 
     and a $\PP^{s}$-bundle over $ P_\PP$. 
     As the morphism $\pi$ is obtained by base changing this family, this establishes the lemma.
    \end{proof}

  As a second step towards proving Claim \ref{cl}, let us identify the generically defined cycles on $F\times Z$:
   
   \begin{lemma}\label{gd} There is equality
    \[    \begin{split} GDA_{B^\circ}^\ast(F\times Z) = \Bigl\langle  (p_{F})^\ast GDA_{B^\circ}^\ast(F)  ,
                  (p_Z)^\ast GDA_{B^\circ}^\ast(Z)\Bigr\rangle&\\
                   +       \Bigl\langle  (p_{F})^\ast GDA_{B^\circ}^\ast(F)  ,
                  (p_Z)^\ast GDA_{B^\circ}^\ast(Z)&\Bigr\rangle \cdot     P \  .\\
                  \end{split}\]
                                  \end{lemma}
                  
   \begin{proof} 
   We have just seen (Lemma \ref{strat}) that 
     \[  \pi\colon\ \  \bar{\FF}\times_{\bar{B}} \bar{\Zz}\ \to\    \Gr(2,6)\times \PP^5\   \  \]
     is a stratified projective bundle, with strata $P_\PP$ and $ \Gr(2,6)\times \PP^5$. Applying \cite[Proposition 5.2]{FLV} to this set-up, we find there is equality
     \begin{equation}\label{sum}  \begin{split}  GDA^\ast_{B^\circ}(F\times Z)= \ima\Bigl( A^\ast( \Gr(2,6)\times \PP^5)\to A^\ast(F\times Z)\Bigr)&\\ +    \iota_\ast \ima\Bigl( A^\ast(P_\PP) \to A^\ast( P)\Bigr)\ ,\\
     \end{split}\end{equation}
     where $\iota\colon  P\hookrightarrow F\times Z$ is the inclusion morphism. The homogeneous varieties $\Gr(2,6)$ and $\PP^5$ have trivial Chow groups and so $A^\ast( \Gr(2,6)\times \PP^5)$ is naturally isomorphic to $A^\ast(\Gr(2,6))\otimes A^\ast(\PP^5)$.
     This means that the first summand of \eqref{sum} can be rewritten as
     \begin{equation}\label{firstsum} \ima\Bigl( A^\ast( \Gr(2,6)\times \PP^5)\to A^\ast(F\times Y)\Bigr)=\Bigl\langle       (p_{F})^\ast GDA_{B^\circ}^\ast(F)  ,
                  (p_Z)^\ast GDA_{B^\circ}^\ast(Z) \Bigr\rangle\   .\end{equation}
     
      As for the second summand of \eqref{sum}, we make the following observation:
      
      \begin{sublemma}\label{leray} The restriction map
        \[  A^\ast\bigl(\Gr(2,6)\times\PP^5\bigr)\ \twoheadrightarrow \ A^\ast(P_\PP)\ \]       
        is surjective.
              \end{sublemma}
              
              \begin{proof} The universal line $P_\PP$  is a $\PP^1$-bundle over $\Gr(2,6)$ with $p^\ast(h)$ relatively ample (where $h\in A^1(\PP^5)$ is ample, and $p\colon P_\PP\to \PP^5$ is the morphism induced by projection). The sublemma thus follows from the projective bundle formula.              
              \end{proof}
        
   Using the surjectivity of Sublemma \ref{leray}, plus the equality \eqref{firstsum}, one reduces \eqref{sum} to the equality of Lemma \ref{gd}. This ends the proof.    
          \end{proof}

As a next step towards proving Claim \ref{cl}, let us make a further simplification to the equality of Lemma \ref{gd}:

\begin{lemma}\label{gd2} There is equality
    \[ GDA_{B^\circ}^\ast(F\times Z) = \Bigl\langle  (p_{F})^\ast GDA_{B^\circ}^\ast(F)  ,
                  (p_Z)^\ast GDA_{B^\circ}^\ast(Z) \Bigr\rangle   \oplus \QQ[P]  \oplus \QQ [P\cdot (p_F)^\ast(h_F)]   \ .\]
                                   \end{lemma}

  \begin{proof} This is Lemma \ref{gd} in combination with the following two sublemmas:
                      
  \begin{sublemma}\label{lemma1}
                 \[ (P)\cdot  (p_Z)^\ast GDA_{B^\circ}^\ast(Z)\ \ \subset\  \Bigl\langle  (p_{F})^\ast GDA_{B^\circ}^\ast(F)  ,
                  (p_Z)^\ast GDA^\ast_{B^\circ}(Z) \Bigr\rangle\ .\]
                               \end{sublemma}
                  
     \begin{proof} We know (Proposition \ref{Fr}) that $GDA_{B^\circ}^\ast(Z)=\langle h\rangle$.
     Also, we know that there exist $a_j\in\QQ$ such that
    \begin{equation}\label{req}
        \Delta_Z\cdot (p_2)^\ast(h) = \sum_j a_j\, h^j\times h^{4-j}\ \ \ \hbox{in}\ A^4(Z\times Z)\ \end{equation}
        \cite[Corollary 3.8]{2q}.
 
  Let $p$ and $q$ denote the projections from $P$ to $F$ resp. $Z$.
  Using equality \eqref{req}, we find that
    \[ \begin{split} P\cdot (p_Z)^\ast(h) &=  (p\times\Delta_Z)_\ast (q\times\Delta_Z)^\ast    \Bigl(\Delta_Z\cdot (p_2)^\ast(h)\Bigr)\\
                         &=   (p\times\Delta_Z)_\ast (q\times\Delta_Z)^\ast  \Bigl(      \sum_j a_j\, h^j\times h^{4-j}\Bigr)\\
                         &= \sum_j a_j\, (p_F)^\ast P^\ast(h^j)\cdot (p_Z)^\ast (h^{4-j})\ \ \ \ \hbox{in}\ A^3(F\times Z)\ .\\
                         \end{split}\]
     Since $P$ and $h$ are generically defined, we have that $P^\ast(h^j)\in GDA_{B^\circ}^\ast(F)$, and so we get
     \[  P\cdot (p_Z)^\ast(h)    \ \ \in \  \Bigl\langle  (p_F)^\ast GDA_{B^\circ}^\ast(F) , 
                  (p_Z)^\ast GDA_{B^\circ}^\ast(Z) \Bigr\rangle\ .\]
    It follows that likewise
     \begin{equation}\label{eqlem} P\cdot (p_Z)^\ast(h^i)    \ \ \in \  \Bigl\langle  (p_F)^\ast GDA_{B^\circ}^\ast(F) , 
                  (p_Z)^\ast GDA_{B^\circ}^\ast(Z) \Bigr\rangle\ \ \ \forall i\ ,\end{equation}
                  which proves the sublemma.
       \end{proof}             
                  
   \begin{sublemma}\label{lemma2}
    \[  (P)\cdot (p_{F})^\ast GDA_{B^\circ}^2(F) \ \ \subset\   \Bigl\langle  (p_{F})^\ast GDA_{B^\circ}^\ast(F)  ,
                  (p_Z)^\ast GDA_{B^\circ}^\ast(Z) \Bigr\rangle    \ .\]      
                   \end{sublemma}
                   
     \begin{proof} In this proof, let us drop the $B^\circ$ subscript, since all generically defined cycles are with respect to $B^\circ$. Since $\FF\to B$ has the Franchetta property, $GDA^1(F)$ and $GDA^2(F)$ are $1$-dimensional, generated by $h_F$ resp. $h_F^2$. It is readily checked (cf. \cite[Sublemma 3.14]{2q}) that $GDA^2(F)$ is also generated by     
     $c_2(Q)\vert_F$, where $Q$ is the universal quotient bundle on $\Gr(2,6)$. To prove the lemma, we thus need to check that
       \begin{equation}\label{inth} P\cdot (p_F)^\ast ( c_2(Q)\vert_F)\ \ \in\   \Bigl\langle  (p_{F})^\ast GDA^\ast(F)  ,
                  (p_Z)^\ast GDA^\ast(Z) \Bigr\rangle     
\ .\end{equation}
     
     The morphism $p\colon P\to F$ being a $\PP^1$-bundle (with $q^\ast(h)$ being relatively ample), we find that
       \[   p^\ast (c_2(Q)\vert_F) = -   q^\ast(h^2) - q^\ast(h) p^\ast (c_1(Q)\vert_F)  \ \ \ \hbox{in}\ A^2(P)\ .\]
       Pushing forward under the closed immersion $P\hookrightarrow F\times Z$, this implies that
       \[  P\cdot (p_F)^\ast ( c_2(Q)\vert_F) =  - P\cdot (p_Z)^\ast (h^2) - P\cdot (p_Z)^\ast(h)\cdot (p_F)^\ast(h_F)\ \ \ \hbox{in}\ A^4(F\times Z)\ .\]
      Using equation \eqref{eqlem}, we see that the right-hand side is decomposable, and so \eqref{inth} is proven. 
     
       (An alternative argument is as follows: up to a finite base change, we may assume $\FF\to B$ is an abelian scheme, i.e. there is a zero section. Let $o\in F$ denote the origin and let $C_o\subset Z$ denote the line corresponding to $o\in F$.
    There is equality
    \[  P\cdot (p_F)^\ast(o) =  (p_F)^\ast(o)\cdot (p_Z)^\ast(C_o)\ \ \ \hbox{in}\ A^4(F\times Z)\ .\]
    Moreover the class $o\in A^2(F)$ is generically defined, and hence is a generator for $GDA^2_B(F)=A^2_{(2)}(F)\cong\QQ$. Also, the class $C_o\in A^2(Z)$
    is generically defined (it is the fiberwise restriction of the relative cycle 
      \[ (p_\Zz)_\ast\bigl( (p_\FF)^\ast(o_\FF)\cdot \PPP\bigr)\ \ \in\ A^2(\Zz)\ ,\] 
      where $o_\FF\in A^2(\FF)$ is the
    zero-section and $\PPP\subset \FF\times_B \Zz$ is the relative universal line). The sublemma is proven.)
    \end{proof}

Combining Sublemmas \ref{lemma1} and \ref{lemma2}, one obtains a proof of Lemma \ref{gd2}.
       \end{proof}
       
     We are now in position to wrap up the proof of Claim \ref{cl}: thanks to Lemma \ref{gd2}, combined with the K\"unneth decomposition in cohomology, the Franchetta property in codimension $\ge 4$ for $\FF\times_{B^\circ}\Zz\to B^\circ$ reduces to the Franchetta property for $\FF\to B^\circ$ and that for $\Zz\to B^\circ$. The second follows from Proposition \ref{Fr2Z}; the first follows from Proposition \ref{Fr2Z} combined with the generically defined isomorphism of the $Z$--$F(Z)$ relation (Proposition \ref{ZF}).
     
For the Franchetta property in codimension 3, one needs to check that
  \begin{equation}\label{cl3}  GDA^3_{B^\circ}(F\times Z)= \dec^3(F\times Z)  \oplus \QQ [P\cdot (p_F)^\ast(h_F)] \ \to\ H^6(F\times Z,\QQ)     \end{equation}
 is injective, where we have used the shorthand
                  \[ \dec^j(F\times Z):=   \Bigl\langle  (p_F)^\ast GDA_{B^\circ}^\ast(F), \,   (p_Z)^\ast GDA_{B^\circ}^\ast(Z) \Bigr\rangle \cap\  A^j(F\times Z)  \]
                  for the {\em decomposable\/} cycles. However, the class $P\cdot (p_F)^\ast(h_F)$ is linearly independent from the decomposable cycles $\dec^3(F\times Z)$ in cohomology, because
            \[  (P\cdot (p_F)^\ast(h_F))_\ast\colon\ \ H^1(F,\QQ)\ \to\ H^3(F,\QQ)\ \to\ H^3(Z,\QQ) \]
       is an isomorphism, while $D_\ast$ is zero on $H^1(F,\QQ)$ for any decomposable cycle $D$. The injectivity of \eqref{cl3} thus reduces to the Franchetta property for $F$ and for $Z$.    
       
       The argument in codimension 2 is similar: one needs to check that
       \begin{equation}\label{cl2}       GDA^2(F\times Z)= \dec^2(F\times Z)  \oplus \QQ [P] \ \to\ H^4(F\times Z,\QQ)    \end{equation}
       is injective. However, the class $P$ is linearly independent from $\dec^2(F\times Z)$ in cohomology, because of the isomorphism
       \[ P_\ast\colon\ \  H^3(F,\QQ)\ \xrightarrow{\cong}\ H^3(Z,\QQ)   \ .\]
       The injectivity of \eqref{cl2} thus reduces to the Franchetta property for $F$ and for $Z$.                                            
                      \end{proof}

    \subsection{Franchetta for $Y^{(2)}\times Y$}
  
  \begin{proposition}\label{Fr3} Let $\YY\to B$ be as above. The family
    \[   \bigl(\YY\times_B \YY\times_B \YY\bigr)/\Sy_2\ \to\ B \]
    (where $\Sy_2$ acts by exchanging the first 2 factors) has the Franchetta property.
  \end{proposition}
  
  \begin{proof} In view of the spread lemma \cite[Lemma 3.2]{Vo}, it suffices to prove the Franchetta property over $B^\circ$.  
  As we have seen, there is an isomorphism of motives
    \[  h(Y)\cong h(Z) \ \ \hbox{in}\ \MM_{\rm rat}\ \]
    which is generically defined with respect to $B^\circ$ (cf. Lemma \ref{ZY} above). This isomorphism of motives induces isomorphisms of Chow groups
    \[  GDA_{B^\circ}^\ast(Y^{(2)}\times Y)\ \xrightarrow{\cong}\ GDA_{B^\circ}^\ast(Z^{(2)}\times Z)\ ,\]
    compatible with cycle class maps. The result thus follows from Proposition \ref{Fr3Z}.
    \end{proof}  

\begin{corollary}\label{FrFY} Let $\YY\to B$ be as above, and let $\FF\to B$ denote the universal family of Fano varieties of conics. The family
  \[ \FF\times_B \YY\ \to\ B \]
  has the Franchetta property.
  \end{corollary}
  
  \begin{proof} This follows from Proposition \ref{Fr3} in view of the generically defined $Y$-$F(Y)$ relation (Proposition \ref{YF}).
  
  \end{proof}
  
  \begin{remark} The Franchetta type result central to this paper (with the purpose of establishing MCK in new cases) is Proposition \ref{Fr3}, which concerns the genus 10 Fano threefold $Y$. However, in order to prove Proposition \ref{Fr3} we were compelled to first prove the analogous result for the index 2 Fano threefold $Z$ (Proposition \ref{Fr3Z}).
  The reason for this detour via $Z$ can be explained as follows: the argument of Proposition \ref{Fr3Z} (roughly speaking: once one has the Franchetta property for $F$ and for $Z$ one also has it for $F\times Z$, as the only ``extra cycles'' come from the universal line $P\subset F\times Z$) does {\em not\/} apply directly to $Y$. Indeed, working with the Fano threefold $Y$ one runs into the double trouble that (1) we do not know whether the correspondence $P$ of Theorem \ref{FY} is the universal conic, (2) 
  even if we knew this, the analogue of Sublemma \ref{leray} is not clear for the universal conic on $Y$.
   \end{remark}

  \section{Main result}
  
  \begin{theorem}\label{main} Let $Y$ be a prime Fano threefold of genus 10. Then $Y$ has a multiplicative Chow--K\"unneth decomposition.
  The induced bigrading on the Chow ring is such that
    \[ A^2_{(0)}(Y) =\QQ[c_2(Y)]\ .\]
  \end{theorem}  
  
\begin{proof} Let $\{\pi^j_Y\}$ be the CK decomposition for $Y$ defined in \eqref{ck}. We observe that this CK decomposition is {\em generically defined\/} with respect to the family $\YY\to B$ (Notation \ref{not}), i.e. it is obtained by restriction from
 ``universal projectors'' $\pi^j_\YY\in A^3(\YY\times_B \YY)$. (This is just because $h$ and $\Delta_Y$ are generically defined.)   
 
 Writing $h^j(Y):=(Y,\pi^j_Y)\in\MM_{\rm rat}$, we have
   \begin{equation}\label{h2j} h^{2j}(Y) \cong \one(-j)\ \ \ \hbox{in}\ \MM_{\rm rat}\ \ \ (j=0,\ldots, 3)\ ,\end{equation}
   i.e. the interesting part of the motive is concentrated in $h^3(Y)$.

 What we need to prove is that this CK decomposition is MCK, i.e.
      \begin{equation}\label{this} \pi_Y^k\circ \Delta_Y^{sm}\circ (\pi_Y^i\times \pi_Y^j)=0\ \ \ \hbox{in}\ A^{6}(Y\times Y\times Y)\ \ \ \hbox{for\ all\ }i+j\not=k\ ,\end{equation}
      or equivalently that
       \[   h^i(Y)\otimes h^j(Y)\ \xrightarrow{\Delta_Y^{sm}}\  h(Y) \]
       coincides with 
       \[ h^i(Y)\otimes h^j(Y)\ \xrightarrow{\Delta_Y^{sm}}\ h(Y)\ \to\ h^{i+j}(Y)  \ \to\ h(Y)\ , \]   
       for all $i,j$.
              
   As a first step, let us assume that we have three integers $(i,j,k)$ and at most one of them is equal to $3$. The cycle in \eqref{this} is generically defined and homologically trivial. The isomorphisms \eqref{h2j} induce an injection
   \[ (\pi^i_Y\times\pi^j_Y\times\pi^k_Y)_\ast A^6(Y\times Y\times Y)\ \hookrightarrow\ A^\ast(Y)\ ,\]
   and send generically defined cycles to generically defined cycles (this is because the isomorphisms \eqref{h2j} are generically defined). As a consequence, the required vanishing \eqref{this} follows from the Franchetta property for $Y$ (Proposition \ref{Fr}).      
   
   In the second step, let us assume that among the three integers $(i,j,k)$, exactly two are equal to $3$. In this case, using the isomorphisms \eqref{h2j} we find an injection
   \[ (\pi^i_Y\times\pi^j_Y\times\pi^k_Y)_\ast A^6(Y\times Y\times Y)\ \hookrightarrow\ A^\ast(Y\times Y)\ ,\] 
respecting the generically defined cycles. As such,  
the required vanishing \eqref{this} follows from the Franchetta property for $Y\times Y$ (Proposition \ref{Fr2}).

In the third and final step, let us treat the case $i=j=k=3$. For this case, we observe that
   \[   \pi_Y^3\circ \Delta_Y^{sm}\circ (\pi_Y^3\times \pi_Y^3)= (\pi^3_Y\times \pi^3_Y\times \pi^3_Y)_\ast (\Delta^{sm}_Y)\ \ \in\  GDA^6_B (Y^{(2)}\times Y)\cap A^6_{hom} (Y^{(2)}\times Y)\ .\]
   The required vanishing \eqref{this} for $i=j=k=3$ thus follows from the Franchetta property for $Y^{(2)}\times Y$ (Proposition \ref{Fr3}).

Finally, let us prove that $A^2_{(0)}(Y)=\QQ[c_2(Y)]$. We remark that $c_2(Y)\in A^2_{(0)}(Y)$ because
   \[   (\pi^i_Y)_\ast c_2(Y) = (\pi^i_\YY)_\ast c_2(T_{\YY/B})\vert_Y = 0\ \ \ \hbox{in}\  A^2(Y)\ \ \  \hbox{for\ all\ } i\not= 4\ ,\]
   as follows from the Franchetta property for $\YY\to B$ (Proposition \ref{Fr}). One readily checks that $c_2(Y)$ is non-zero (e.g. one can take a smooth anticanonical section $S\subset Y$; if $c_2(Y)$ were zero then by adjunction also $c_2(S)=0$, which is absurd since $S$ is a K3 surface). Since $A^2_{(0)}(Y)$
   injects into $H^4(Y,\QQ)\cong\QQ$, it follows that $A^2_{(0)}(Y)=\QQ[c_2(Y)]$.   
    \end{proof}

 \section{Compatibilities}
 
 In this short section, we show that the MCK decomposition we have constructed for $Y$ is compatible with the one on the Fano surface $F(Y)$ (Proposition \ref{compat}), and with the one on the associated index 2 Fano threefold $Z$ (Proposition \ref{compat2}).

  \begin{proposition}\label{compat} Let $Y$ be a prime Fano threefold of genus 10, let $F:=F(Y)$ be the Fano surface of conics and let $P\subset F\times Y$ be the universal conic. Then 
    \[   P_\ast A^i_{(j)}(F)\ \subset\ A^{i-1}_{(j)}(Y)\ ,\ \ \   P^\ast A^i_{(j)}(Y)\ \subset\ A^i_{(j)}(F) \ ,\]
    where $A^\ast_{(\ast)}(Y)$ is the bigrading induced by the MCK decomposition of Theorem \ref{main}, and $A^\ast_{(\ast)}(F)$ is the Beauville bigrading for the abelian surface $F$.    
  \end{proposition}
  
  \begin{proof} 
  Let $\pi^\ast_Y$ and $\pi^\ast_F$ denote the MCK decomposition of Theorem \ref{main} resp. the MCK decomposition of the abelian variety $F$. We will prove that the correspondence $P$ is of pure degree 0, in the sense of Shen--Vial \cite[Definition 1.2]{SV2}, i.e.
    \begin{equation}\label{van} (\pi_F^i\times \pi_Y^j)_\ast P = 0\ \ \ \forall\ i+j\not= 6\ .\end{equation}
 This implies the statement in view of \cite[Proposition 1.6]{SV2}.
 
 To prove the vanishing \eqref{van}, we observe that the cycle in \eqref{van} is generically defined (with respect to the base $B$) and homologically trivial. The vanishing thus follows from the Franchetta property for $F\times Y$ (Corollary \ref{FrFY}).
   \end{proof}

  \begin{proposition}\label{compat2} Let $Y$ be a prime Fano threefold of genus 10, and let $Z$ be the index 2 Fano threefold associated to $Y$ (in the sense that $F_1(Z)\cong F(Y)$, cf. Theorem \ref{FY}). Then there are isomorphisms
   \[  A^i_{(j)}(Y)\cong A^i_{(j)}(Z)\ ,\] 
   where $A^\ast_{(\ast)}(Y)$ is the bigrading induced by the MCK decomposition of Theorem \ref{main}, and $A^\ast_{(\ast)}(Z)$ is the bigrading constructed in \cite{2q}.
   \end{proposition}
   
   \begin{proof} Let $\pi^\ast_Y$ and $\pi^\ast_Z$ denote the MCK decomposition of Theorem \ref{main} resp. the MCK decomposition of \cite{2q}; both are generically defined with respect to $B^\circ$.  
   As in Proposition \ref{compat}, it suffices to prove that the correspondence $ $ of Lemma \ref{ZY} is of pure degree 0. This follows from the Franchetta property for $Y\times Z$, which is equivalent to the Franchetta property for $Y\times Y$ in view of Lemma \ref{ZY}.
     \end{proof}

 \section{The tautological ring}
 
 \begin{corollary}\label{cor1} Let $Y$ be a prime Fano threefold of genus 10, and let $m\in\NN$. Let
  \[ R^\ast(Y^m):=\bigl\langle (p_i)^\ast(h), (p_{ij})^\ast(\Delta_Y)\bigr\rangle\ \subset\ \ \ A^\ast(Y^m)   \]
  be the $\QQ$-subalgebra generated by pullbacks of the polarization $h\in A^1(Y)$ and pullbacks of the diagonal $\Delta_Y\in A^3(Y\times Y)$. (Here $p_i$ and $p_{ij}$ denote the various projections from $Y^m$ to $Y$ resp. to $Y\times Y$).
  The cycle class map induces injections
   \[ R^\ast(Y^m)\ \hookrightarrow\ H^\ast(Y^m,\QQ)\ \ \ \hbox{for\ all\ }m\in\NN\ .\]
   \end{corollary}

\begin{proof} This is inspired by the analogous result for cubic hypersurfaces \cite[Section 2.3]{FLV3}, which in turn is inspired by analogous results for hyperelliptic curves \cite{Ta2}, \cite{Ta} (cf. Remark \ref{tava} below) and for K3 surfaces \cite{Yin}.

As in \cite[Section 2.3]{FLV3}, let us write $o:={1\over 18} h^3\in A^3(Y)$, and
  \[ \tau:= \Delta_Y - {1\over 18}\, \sum_{j=0}^3  h^j\times h^{3-j}\ \ \in\ A^3(Y\times Y) \]
  (this cycle $\tau$ is nothing but the projector on the motive $h^3(Y)$ considered above).
Moreover, let us write 
  \[ \begin{split}   o_i&:= (p_i)^\ast(o)\ \ \in\ A^3(Y^m)\ ,\\
                        h_i&:=(p_i)^\ast(h)\ \ \in \ A^1(Y^m)\ ,\\
                         \tau_{i,j}&:=(p_{ij})^\ast(\tau)\ \ \in\ A^3(Y^m)\ .\\
                         \end{split}\]
We define the $\QQ$-subalgebra
  \[ \bar{R}^\ast(Y^m):=\langle o_i, h_i, \tau_{ij}\rangle\ \ \ \subset\ H^\ast(Y^m,\QQ) \]
  (where $i$ ranges over $1\le i\le m$, and $1\le i<j\le m$). One can prove (just as \cite[Lemma 2.11]{FLV3} and \cite[Lemma 2.3]{Yin}) that the $\QQ$-algebra $ \bar{R}^\ast(Y^m)$
  is isomorphic to the free graded $\QQ$-algebra generated by $o_i,h_i,\tau_{i,j}$, modulo the following relations:
    \begin{equation}\label{E:X'}
			o_i\cdot o_i = 0, \quad h_i \cdot o_i = 0,  \quad 
			h_i^3 =18\,o_i\,;
			\end{equation}
			\begin{equation}\label{E:X2'}
			\tau_{i,j} \cdot o_i = 0 ,\quad \tau_{i,j} \cdot h_i = 0, \quad \tau_{i,j} \cdot \tau_{i,j} = 4\, o_i\cdot o_j
			\,;
			\end{equation}
			\begin{equation}\label{E:X3'}
			\tau_{i,j} \cdot \tau_{i,k} = \tau_{j,k} \cdot o_i\,;
			\end{equation}
			\begin{equation}\label{E:X4'}
			\sum_{\sigma \in \mathfrak{S}_{6}}  \prod_{i=1}^{6} \tau_{\sigma(2i-1), \sigma(2i)} = 0\,. 
			\end{equation}

To prove Corollary \ref{cor1}, we need to check that these relations are also verified modulo rational equivalence.
The relations \eqref{E:X'} take place in $R^\ast(Y)$ and so they follow from the Franchetta property for $Y$ (Proposition \ref{Fr}). 
The relations \eqref{E:X2'} take place in $R^\ast(Y^2)$. The first and the last relations are trivially verified, because $Y$ being Fano one has $A^6(Y^2)=\QQ$. As for the second relation of \eqref{E:X2'}, this follows from the Franchetta property for $Y\times Y$ (Proposition \ref{Fr}(\rom2)). (Alternatively, one can deduce the second relation from the MCK decomposition: the product $\tau_{} \cdot h_i$ lies in $A^4_{(0)}(Y^2)$, and it is readily checked that $A^4_{(0)}(Y^2)$ injects into $H^8(Y^2,\QQ)$.)
   
   Relation \eqref{E:X3'} takes place in $R^\ast(Y^3)$ and follows from the MCK relation. Indeed, we have
   \[  \Delta_Y^{sm}\circ (\pi^3_Y\times\pi^3_Y)=   \pi^6_Y\circ \Delta_Y^{sm}\circ (\pi^3_Y\times\pi^3_Y)  \ \ \ \hbox{in}\ A^6(Y^3)\ ,\]
   which (using Lieberman's lemma) translates into
   \[ (\pi^3_Y\times \pi^3_Y\times\Delta_Y)_\ast    \Delta_Y^{sm}  =   ( \pi^3_Y\times \pi^3_Y\times\pi^6_Y)_\ast \Delta_Y^{sm}                                                         \ \ \ \hbox{in}\ A^6(Y^3)\ ,\]
   which means that
   \[  \tau_{1,3}\cdot \tau_{2,3}= \tau_{1,2}\cdot o_3\ \ \ \hbox{in}\ A^6(Y^3)\ .\]
   
%
 
  Finally, relation \eqref{E:X4'}, which takes place in $R^\ast(Y^{12})$, 
  expresses the fact that 
     \[ \sym^{6}  H^3(Y,\QQ)=0\ ,\]
     where $H^3(Y,\QQ)$ is viewed as a super vector space.
 To check that this relation is also verified modulo rational equivalence, 
 we observe that relation \eqref{E:X4'} involves a cycle  
 contained in
   \[ A^\ast(\sym^{6} (h^3(Y))\ .\]
 But we have vanishing of the Chow motive
   \[ \sym^{6} h^3(Y)=0\ \ \hbox{in}\ \MM_{\rm rat}\ ,\] 
 because $\dim H^3(Y,\QQ)=4$ and $h^3(Y)$ is oddly
  finite-dimensional (all Fano threefolds have finite-dimensional motive \cite[Theorem 4]{43}).  
  This establishes relation \eqref{E:X4'}, modulo rational equivalence, and ends the proof.

%
%
%
 \end{proof}

\begin{remark}\label{tava} Given a curve $C$ and an integer $m\in\NN$, one can define the {\em tautological ring\/}
  \[ R^\ast(C^m):= \bigl\langle  (p_i)^\ast(K_C),(p_{ij})^\ast(\Delta_C)\bigr\rangle\ \ \ \subset\ A^\ast(C^m) \]
  (where $p_i, p_{ij}$ denote the various projections from $C^m$ to $C$ resp. $C\times C$).
  Tavakol has proven \cite[Corollary 6.4]{Ta} that if $C$ is a hyperelliptic curve, the cycle class map induces injections
    \[  R^\ast(C^m)\ \hookrightarrow\ H^\ast(C^m,\QQ)\ \ \ \hbox{for\ all\ }m\in\NN\ .\]
   On the other hand, there are many (non hyperelliptic) curves for which the tautological ring $R^\ast(C^3)$ does {\em not\/} inject into cohomology (this is related to the non-vanishing of the Ceresa cycle, cf. \cite[Remark 4.2]{Ta} and also \cite[Example 2.3 and Remark 2.4]{FLV2}). 
   
Corollary \ref{cor1} shows that genus 10 Fano threefolds behave similarly to hyperelliptic curves. It would be interesting to understand what happens for other  Fano threefolds: is Corollary \ref{cor1} true for all of them or not ? (This is related to Question \ref{ques}.)
\end{remark}

\section{Question}
\label{s: q}

\begin{question} Let $Y$ be a prime Fano threefold of genus $10$, and let $S\subset Y$ be a smooth anticanonical divisor. Then $S$ is a K3 surface (and the general genus 10 K3 surface arises in this way). Is it true that
  \[ \begin{split} \ima\bigl(  A^1(S)\ \to\ A^2(Y)\bigr) &= \QQ[h^2]\ ,\\
                              \ima\bigl(  A^2(Y)\ \to\ A^2(S)\bigr) &= \QQ[o_S] \ \ ??\\
  \end{split}\]
 (To prove this, it would suffice to prove that the graph of the inclusion morphism $S\hookrightarrow Y$ is in $A^3_{(0)}(S\times Y)$, with respect to the product MCK decomposition. I have not been able to settle this.)
 
 This question also makes sense for other Fano threefolds and their anticanonical divisors (e.g., the question is interesting for cubic threefolds). 
  \end{question}

 \vskip1cm
\begin{nonumberingt} Thanks to Laurent Manivel for kindly answering my questions and alerting me to \cite{Man}. Thanks to Kai and Len for being great experts on Fort Boyard.
\end{nonumberingt}

\end{document}